\documentclass[12pt]{amsart}
\usepackage{epsfig,psfrag,verbatim,amssymb,amsfonts}
\newtheorem{theorem}{Theorem}
\newtheorem{lemma}[theorem]{Lemma}
\newtheorem{prop}[theorem]{Proposition}

\newcommand{\field}[1]{\mathbb{#1}}
\newcommand{\R}{\field{R}}

\newcommand{\F}{{\mathcal F}}

\newcommand{\si}{\sigma}
\newcommand{\al}{\alpha}
\newcommand{\ga}{\gamma}

\renewcommand{\phi}{\varphi}

 \newcommand{\won}{{\boldsymbol 1}}
\newcommand{\foot}{}

\newcounter{constante}
\newcommand{\con}[1]{
                    \immediate\write 1{\noexpand\newlabel{#1}{{\theconstante}{\theconstante}}}
                    c_{\theconstante}
                    \stepcounter{constante}
                   }
\newcommand{\abel}[1]{}

\begin{document}

\setcounter{page}{1}

\title[Random walk avoiding its past convex hull]
{On the speed of a planar random walk avoiding its past convex hull}
\thanks{\textit{2000 Mathematics Subject Classification.} 60K35.
Secondary 60G50, 52A22.}
\thanks{\textit{Key words:}\quad convex hull, large deviations,
law of large numbers,
 random walk,
 self avoiding, speed}

\maketitle

\begin{center}
{\sc By Martin P.W.\ Zerner}

\end{center}\vspace*{5mm}

\begin{center}\begin{quote}
{{\sc Abstract}. \small We consider a  random walk in $\R^2$
 which takes steps uniformly distributed on the unit circle centered around the walker's
current position but avoids the convex hull of its past positions.
This model has been introduced by Angel, Benjamini and
Vir\'{a}g. We  show a large deviation estimate for the distance of the walker
from the origin, which implies that the walker has positive lim inf speed. 
}\vspace*{5mm}
\end{quote}
\end{center}

\section{Introduction}
\begin{figure}[t]
\epsfig{file=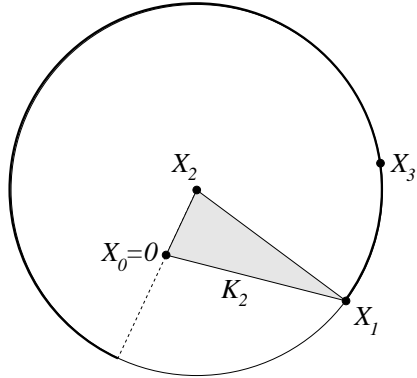,height=5cm,angle=0}
\caption{\foot Three steps of the walk.
$X_3$ is uniformly distributed on the bold arc of the circle with radius 1, centered in $X_2$.
}\label{def}\vspace*{-3mm}
\end{figure}

Angel, Benjamini and Vir\'{a}g introduced in \cite{ABV} the following model
of a random walk $(X_n)_{n\geq 0}$ in $\R^2$, which they  called the \textit{rancher}.
The walker starts at the origin $X_0=0$.
Suppose it  
has already taken $n$ steps ($n\geq 0$) and is currently
at $X_n$. 
Then its next position $X_{n+1}$ is uniformly distributed on the unit circle
centered around $X_n$ but conditioned so that the straight line segment
$\overline{X_n,X_{n+1}}$ from $X_n$ to $X_{n+1}$ does not intersect
the interior $K_n^o$ of the convex hull $K_n$ of 
the past positions $\{X_0,X_1,\ldots,X_n\}$, see Figure \ref{def}.
\abel{def}
Note that $(X_n)_{n\geq 0}$ is not Markovian since in general one needs to know the whole
history of the process in order to determine the transition probabilities
for the next step. This makes this model difficult to analyse, a property it shares with
many other self-interacting processes, see \cite{ABV} and also \cite{BW} for references.

To the best of our knowledge, the only major rigorous result which has been proved
so far for this model, see \cite[Theorem 1]{ABV},
is  that the walk has positive lim sup
speed, i.e.\  there is a constant $c>0$ such that $P$-a.s.\ 
$\limsup \|X_n\|/n
>c$ as $n\to\infty$.  Here  $(\Omega, \F, P)$
is the underlying probability space. 

The purpose of the present paper is to improve this result by showing the 
following.
\begin{theorem}\label{kalt}
There is a constant $\con{one}>0$ such that
\begin{equation}\label{ldp}
\limsup_{n\to\infty}\frac{1}{n}\log P\left[\|X_n\|\leq  c_{\ref{one}}n\right]<0
\end{equation}
and consequently,
\begin{equation}\label{part}
\liminf_{n\to\infty}\frac{\|X_n\|}{n}\geq c_{\ref{one}}\qquad\mbox{$P$-a.s..}
\end{equation}
\end{theorem}\abel{kalt}\abel{one}\abel{part}
In particular, (\ref{part}) proves \cite[Conjecture 4]{ABV}.
We expect but were not able to prove that the speed $\lim\|X_n\|/n$
exists and is $P$-a.s.\ constant, as conjectured in \cite[Conjecture 5]{ABV}. 
For more conjectures regarding convergence of $X_n/\|X_n\|$ and transversal fluctuations
 of trajectories, see \cite{ABV}.

Let us now describe how the present article is organized. The next section introduces
general notation. In Section 3 we introduce some sub- and supermartingales, which 
enable us in Section 4  to
bound exponential moments of the time it takes the diameter
of the convex hull $K_n$ to increase. From this we deduce in Section 5 estimates for the diameter of $K_n$
 similar to the ones claimed  in Theorem \ref{kalt} for $\|X_n\|$
and show how this implies Theorem \ref{kalt}.

%%%%%%%%%%%%%%%%%%%%%%%%%%%%%%%%%%%%%%%%%%%%%%%%%%%%%%%%%%%%%%%%%%%%%%

\section{Notation} 
\begin{figure}[t]
\epsfig{file=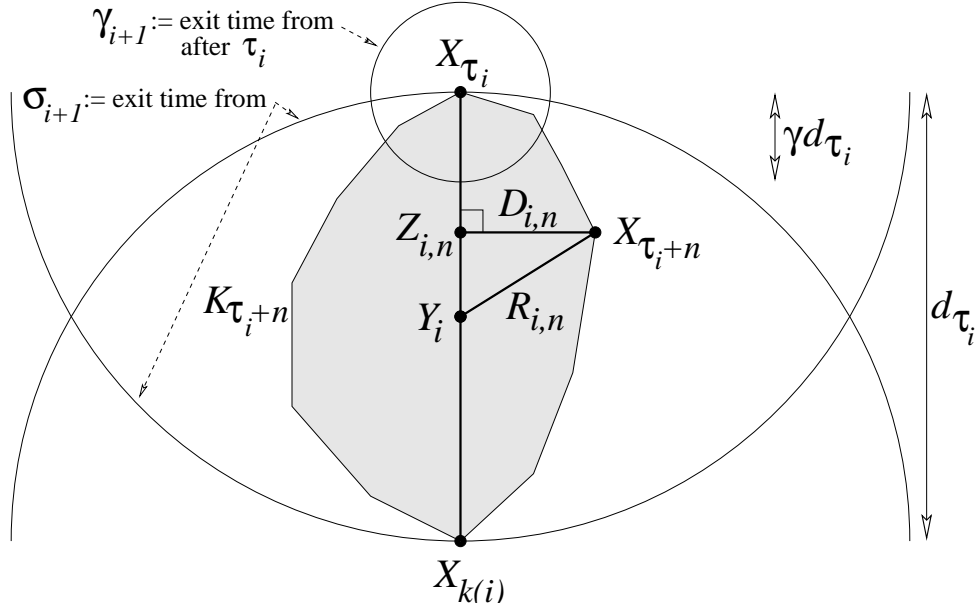,height=8cm,angle=0}
\caption{\foot General notation
}\label{weih2}\vspace*{-2mm}
\end{figure}

\begin{figure}[t]
\epsfig{file=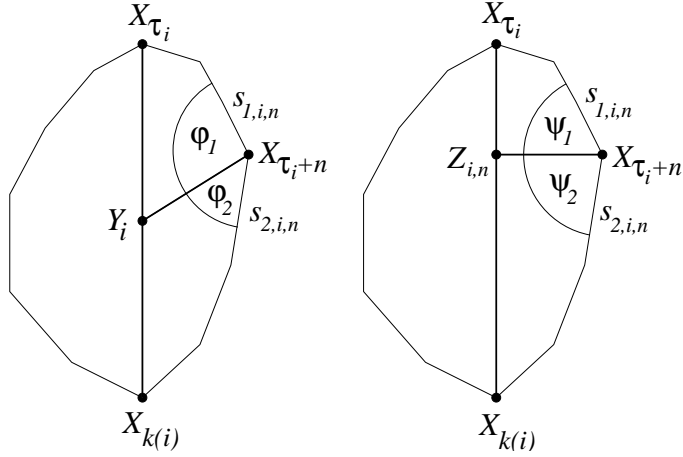,height=6cm,angle=0}
\caption{\foot General notation
}\label{weih}\vspace*{-4mm}
\end{figure}

 We denote by 
$d_n$
the diameter of $K_n$.
Since $(K_n)_n$ is an increasing sequence of sets, $(d_n)_n$ is non-decreasing.
The ladder times $\tau_i$ at which the process $(d_n)_{n\geq 0}$  strictly increases 
are defined
recursively by
\[\tau_0:=0\quad \mbox{and}\quad
\tau_{i+1}:=\inf\{n>\tau_i\ :\ d_n>d_{\tau_i}\}\quad(\leq \infty)\qquad(i\geq 0).
\]
Note that $\tau_1=1$ and that the $\tau_i$'s    are stopping times with respect to the
canonical filtration $(\F_n)_{n\geq 0}$ generated by $(X_n)_{n\geq 0}$.
Since the diameter of a bounded convex set is the distance between two of its extremal points
there is for all $i\geq 1$ with $\tau_i<\infty$ a ($P$-a.s.\ unique) $0\leq k(i)< \tau_i$ such that
$d_{\tau_i}=\|X_{\tau_i}-X_{k(i)}\|$, see Figure \ref{weih2}.
\abel{weih2}
For $x\in\R^2$ and $r>0$ we denote by $B(x,r)$ the closed disk with center $x$ and radius $r$.  
If $\tau_i<\infty$ then
\[
\si_{i+1}:=\inf\{n\geq 0\mid X_n\notin B(X_{\tau_i}, d_{\tau_i})\cap
B(X_{k(i)}, d_{\tau_i})\}
\]
is the exit time of the walk 
from the large  lens shaped region shown in Figure \ref{weih2}, which we shall
refer to as the \textit{lens created at time $\tau_i$}. Observe that $K_{\tau_i}$ is contained
in the lens created at time $\tau_i$. Moreover,
\begin{equation}\label{stau}
\tau_{i+1}\leq \si_{i+1}
\end{equation}\abel{stau}
since if $\si_{i+1}<\infty$, $X_{\si_{i+1}}$ has a  distance from either
$X_{\tau_i}$ or $X_{k(i)}$ greater than  
 $d_{\tau_i}$.
The point
\begin{equation}\label{above}
Y_i:=\frac{X_{\tau_i}+X_{k(i)}}{2}\qquad (i\geq 1)
\end{equation}\abel{above}
will serve as the ``center" of $K_{\tau_i}$ and 
\[R_{i,n}:=\left\|X_{\tau_i+n}-Y_i\right\|\qquad (i\geq 1, n\geq 0)\]
is the  distance of $X_{\tau_i+n}$ from this center.
The orthogonal projection of $X_{\tau_i+n}$ onto the straight line passing through 
$X_{\tau_i}$ and $X_{k(i)}$ will be called $Z_{i,n}\ (i\geq 1, n\geq 0)$. 
The distance of $X_{\tau_i+n}$ from this line 
is denoted by
\[D_{i,n}:=\|X_{\tau_i+n}-Z_{i,n}\|\qquad (i\geq 1, n\geq 0).
\]
For the  following definitions we assume $i,n\geq 1$ and $\tau_i+n<\tau_{i+1}$. 
In particular, due to (\ref{stau}), we assume 
that at time $\tau_i+n$ the walk has not yet left the lens created at time $\tau_i$ .
This implies that
$Z_{i,n}\in \overline{X_{\tau_i},X_{k(i)}}$ and that 
$X_{\tau_i}$ and $X_{k(i)}$ are still boundary points
of $K_{\tau_i+n}$, as shown in Figures \ref{weih2} and \ref{weih}.
Hence if we start in $X_{\tau_i+n}$ and follow the two boundary line segment emanating from $X_{\tau_i+n}$
we will eventually reach $X_{\tau_i}$ and  $X_{k(i)}$. The boundary line segment whose continuation  leads first
 to $K_{\tau_i}$ and then to $X_{k(i)}$ is called $s_{1,i,n}$,  
while the other  line segment starting in $X_{\tau_i+n}$ is denoted by $s_{2,i,n}$, see Figure \ref{weih}.
\abel{weih}
The angle between $s_{j,i,n}$ and 
$\overline{X_{\tau_i+n},Y_{i}}$ 
 is
called $\phi_{j,i,n}\in[0,\pi]$ $(j=1,2)$, see the left part
of Figure \ref{weih}. 
Similarly, the angle between $s_{j,i,n}$ and $\overline{X_{\tau_i+n},Z_{i,n}}$ is denoted by 
$\psi_{j,i,n}\in[0,\pi]$ $(j=1,2)$, see the right part
of Figure \ref{weih}. 
Occasionally, we will dropped the subscripts $i$ and $n$ from $\phi$ and $\psi$.
Since $K_{\tau_i+n}$ is convex, 
\begin{equation}\label{Code}
\phi_1+\phi_2= \psi_1+\psi_2\leq \pi.
\end{equation}\label{code}
Furthermore,  $|\phi_1-\psi_1|$ is
one of the  angles in a right angled
triangle, namely the triangle with vertices $X_{\tau_i+n}, Y_i$ and $Z_{i,n}$. Hence,
\begin{equation}\label{fax}
|\phi_1-\psi_1|= |\phi_2-\psi_2|\leq\pi/2.
\end{equation}\abel{fax}

%%%%%%%%%%%%%%%%%%%%%%%Some sub- and supermartingales%%%%%%%%%%%%%%%%%%%%%%%%%%%%%%%%%%

\section{Some sub- and supermartingales}

\begin{figure}[t]
\epsfig{file=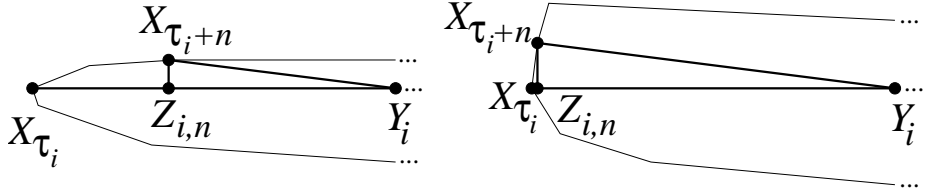,height=2.5cm,angle=0}
\caption{\foot  The expected increment of $R_{i,n}$ is small in the left figure and large
in the right figure. For $D_{i,n}$ it is the other way round. 
}\label{bounded}\vspace*{-4mm}
\end{figure}

The following result shows that for every $i\geq 1$,  both 
$(R_{i,n})_n$ and $(D_{i,n})_n\ (1\leq n<\tau_{i+1}-\tau_i)$ 
are submartingales.
\begin{lemma}\label{neg}
For all $i, n\geq 1$, $P$-a.s.\ on $\{\tau_i+n<\tau_{i+1}\}$,
\begin{eqnarray}
E[R_{i,n+1}- R_{i,n}\ |\ \F_{\tau_i+n}]\label{blitz}
&\geq& \frac{\sin \phi_{1,i,n}+\sin\phi_{2,i,n}}{2\pi}\geq \frac{\sin \phi_{1,i,n}}{2\pi}\geq 0,\\
E[D_{i,n+1}-D_{i,n}\ |\ \F_{\tau_i+n}]\label{nosea}
&\geq& \frac{\sin \psi_{1,i,n}+\sin \psi_{2,i,n}}{2\pi}\ \geq\ \frac{\sin \psi_{1,i,n}}{2\pi}\geq 0
\end{eqnarray}
and
\begin{equation}
E[D_{i,n+1}-D_{i,n}+R_{i,n+1}-R_{i,n}\ |\ \F_{\tau_i+n}]\geq c_{\ref{ekg}}
\label{jodea}
\end{equation}
for some constant $\con{ekg}>0$. \abel{ekg}
\end{lemma}\abel{blitz}\abel{nosea}\abel{neg}\abel{jodea}
%%%%%%%%%%%%%%%%%%%%%Figure boundedaway%%%%%%%%%%%%%%%%%%%%%%%%%%%%%%%%%%%%%%%%%%%%%%%%%%
Figure \ref{bounded}
shows examples in which the  expected increments of   $R_{i,n}$ and $D_{i,n}$
are close to 0, thus explaining, why we are not able to bound 
in (\ref{blitz}) and (\ref{nosea}) these expected increments
individually away from 0.
Note however, that in both situation depicted in Figure \ref{bounded}, if the expected increment
of $R_{i,n}$ or of $D_{i,n}$ is small then the expected increment of the other quantity is large.
This confirms that the expected increments of $R_{i,n}$ and $D_{i,n}$ cannot both be small
at the same time, see (\ref{jodea}).
\abel{bounded}
\begin{proof}[Proof of Lemma \ref{neg}] We
fix $i, n\geq 1$ and drop them as subscripts of $\phi_{j,i,n}$ and $\psi_{j,i,n}$ $(j=1,2)$.
Then the following statements hold  on the event  $\{\tau_i+n<\tau_{i+1}\}$. 
Consider the angle  between $\overline{Y_i,X_{\tau_i+n}}$ and  
$\overline{X_{\tau_i+n},X_{\tau_i+n+1}}$ 
which includes $s_{2,i,n}$. 
This angle
is chosen uniformly at random from the interval $[\phi_2,2\pi-\phi_1]$.
Hence we get by a change of basis argument 
\begin{eqnarray*}
\lefteqn{E[R_{i,n+1}\ |\ \F_{\tau_i+n}]}\\
&=&\frac{1}{(2\pi-\phi_1)-\phi_2}\int_{\phi_2}^{2\pi-\phi_1}
\left\|(R_{i,n},0)-(\cos \phi,\sin\phi)\right\|\ d\phi\\
&\geq& \frac{1}{2\pi-\phi_1-\phi_2}\int_{\phi_2}^{2\pi-\phi_1}|R_{i,n}-\cos \phi|\ d\phi\\
&\geq& R_{i,n}+ \frac{1}{2\pi-\phi_1-\phi_2}\int_{\phi_2}^{2\pi-\phi_1}-\cos \phi\ d\phi\\
&=&  R_{i,n}+ \frac{\sin \phi_1+\sin\phi_2}{2\pi-\phi_1-\phi_2}\ \geq\
 R_{i,n}+ \frac{\sin \phi_1+\sin\phi_2}{2\pi},
\end{eqnarray*}
which shows (\ref{blitz}).
Similarly, (\ref{nosea}) follows from 
\[E[D_{i,n+1}\ |\ \F_{\tau_i+n}
]= \frac{1}{(2\pi-\psi_1)-\psi_2}\int_{\psi_2}^{2\pi-\psi_1}
|D_{i,n}-\cos \psi|\ d\psi.
\]
For the proof of (\ref{jodea})
we assume without loss of generality $\phi_1\leq \pi/2$. Indeed, otherwise
 $\phi_2\leq \pi/2$ because of $\phi_1+\phi_2\leq \pi$, see (\ref{Code}),
and in the following proof  one only has to replace the subscript 1 by the subscript 2 and swap
$X_{\tau_i}$
and $X_{k(i)}$.
By (\ref{blitz}) and 
(\ref{nosea}),
\begin{eqnarray}\label{zhan}
E[D_{i,n+1}-D_{i,n}+R_{i,n+1}-R_{i,n}\ |\ \F_{\tau_i+n}]&\geq& 
\frac{\sin \phi_1+\sin \psi_1}{2\pi}.
\end{eqnarray}\abel{zhan}
We will show that 
the right side of (\ref{zhan}) is
always greater than $c_{\ref{ekg}}:=(4\pi^2)^{-1}$.
Assume that it is less than  $c_{\ref{ekg}}$. Then 
\begin{equation}\label{warum}
\sin \phi_1, \sin \psi_1\leq (2\pi)^{-1}
\end{equation}\abel{warum}
 and hence $\phi_1\leq (\pi/2)\sin\phi_1\leq 1/4$ 
by concavity of $\sin$ on $[0,\pi/2]$.
Similarly, (\ref{warum}) implies that either $\psi_1\leq 1/4$ or $\pi-\psi_1\leq 1/4$.
Due to $|\phi_1-\psi_1|\leq \pi/2$, see (\ref{fax}), the latter case is impossible.
 Therefore, 
\begin{equation}\label{all}
|\phi_1-\psi_1|\leq \max\{|\phi_1|,|\psi_1|\}\leq 1/4.
\end{equation}\abel{all}
The angle $\al\in [0,\pi/2]$ 
between $\overline{Z_{i,n},X_{\tau_i+n}}$ and $\overline{X_{\tau_i+n},X_{\tau_i}}$
is less than or equal to $\psi_1$. Consequently,
\begin{equation}\label{oaks}
\sin\psi_1\ \geq\ \sin\al\ =\ \frac{\|X_{\tau_i}-Z_{i,n}\|}{\|X_{\tau_i}-X_{\tau_i+n}\|}\geq 
\frac{\|X_{\tau_i}-Y_i\|-\|Y_i-Z_{i,n}\|}{d_{\tau_i}}.
\end{equation}\abel{oaks}
However,
\begin{equation}\label{nochmal}
R_{i,n}=\|Y_i-X_{\tau_i+n}\|\stackrel{(\ref{above})}{\leq}
 (\|X_{\tau_i}-X_{\tau_i+n}\|+\|X_{k(i)}-X_{\tau_i+n}\|)/2
\leq d_{\tau_i}.
\end{equation}\abel{nochmal}
Therefore, 
\[
\sin\psi_1\stackrel{(\ref{oaks})}{\geq} \frac{d_{\tau_i}/2}{d_{\tau_i}}-\frac{\|Y_i-Z_{i,n}\|}{\|Y_i-X_{\tau_i+n}\|}
\ =\ \frac{1}{2}-\sin|\phi_1-\psi_1|
\ \stackrel{(\ref{all})}{\geq}\  \frac{1}{2}-\frac{1}{4}=\frac{1}{4},
\]
which contradicts (\ref{warum}).
\end{proof}
%%%%%%%%%%%%%%%%%%%%%%%%%%%%%%%%%%%%%%%%%%%%%%%%%%%%%%%%%%%%%%%%%%%%%%%

We fix the constants
\begin{equation}\label{bega}
\beta:=1+4\pi \sqrt{8}> 30\quad\mbox{and}\quad\gamma:=\frac{1}{2\beta}<\frac{1}{60}.
\end{equation}\abel{bega}
Whenever $\tau_i<\infty$ we denote the 
first exit time  after $\tau_i$ from  $B(X_{\tau_i}, \ga d_{\tau_i})$ by
\[
\ga_{i+1}:=\inf\{n>\tau_i\ :\ \|X_n-X_{\tau_i}\|>\ga  d_{\tau_i} \}\quad (\leq \infty),
\]
see Figure \ref{weih2}.
If $i\geq 1$ and $n\geq 0$  then we  shall call  $n$ \textit{good for} $i$
if $n=0$ or if 
\begin{equation}\label{lean}
\tau_i+n< \tau_{i+1}\wedge\ga_{i+1}\quad\mbox{and}
\quad E[R_{i,n+1}-R_{i,n}\ |\  \F_{\tau_i+n}]\geq \frac{1}{\pi\sqrt{8}}\ \mbox{$P$-a.s..}
\end{equation}\abel{lean}
This means, $n\geq 1$ is good for $i$ if at time $\tau_i+n$ the walker 
has not yet left 
the intersection of the small ball around $X_{\tau_i}$ 
and the lens shown in Figure \ref{weih2} and, roughly speaking,
feels a substantial centrifugal force pushing it away from the center $Y_i$.  
Good times help the walker to leave the lens shortly after $\tau_i$ and  closely to the point
$X_{\tau_i}$. Next we introduce a family of supermartingales,
which will help make this idea more precise.

%%%%%%%%%%%%%%%%%%%%%%%%%%%%%%%%%%%%%%%%%%%%%%%%%%%%%%%%%%%%%%%%%%%%%%%%%%%%%%%%%%%%%%

\begin{lemma}\label{ftau}
There are constants $\con{alll}>0, \con{ofer}>0$ and $1\leq\con{elon}<\infty$ such that $P$-a.s.\  for all $i\geq 1, n\geq 0$,
\begin{equation}
\label{chen}
E\left[M_{i,n}\ |\ \F_{\tau_i}\right]\leq c_{\ref{elon}}\exp(-c_{\ref{ofer}}n),
\end{equation}
where 
\begin{eqnarray}\nonumber
M_{i,n}&:=&\won\{\tau_i+n<\tau_{i+1}\} \\
&&\times \exp\bigg(-c_{\ref{alll}}\bigg(D_{i,n}+\beta  (R_{i,n}-R_{i,0})-4\sum_{j=0}^{n-1} 
\won\{\mbox{\rm $j$ is good for $i$}\}\bigg)\bigg).\label{mart}
\end{eqnarray}\abel{chen} \abel{mart}\abel{ftau}\abel{elon}
\end{lemma}\abel{ftau}\abel{elon}
\begin{proof} 
Firstly, we shall prove that for suitable $c_{\ref{ofer}}>0$,
 \begin{equation}\label{chin}
E\left[M_{i,n+1}\ |\ \F_{\tau_i+n}\right]\leq \exp(-c_{\ref{ofer}})M_{i,n}\qquad
\mbox{$P$-a.s. for all $i\geq 1, n\geq 1$,}
\end{equation}\abel{ofer}\abel{chin}
thus showing that $(M_{i,n})_{n\geq 1}$ is an exponentially fast decreasing submartingale for
each $i\geq 1$. Fix $i\geq 1$ and $n\geq 1$.
We have
\begin{eqnarray*}
E\left[M_{i,n+1}\mid \F_{\tau_i+n}\right]&\leq&M_{i,n}f_{i,n}(c_{\ref{alll}})\quad\mbox{$P$-a.s., where}\\
f_{i,n}(c_{\ref{alll}})&:=&E\left[\exp\left(c_{\ref{alll}} Z_{i,n}\right) 
 \ |\ \F_{\tau_i+n}\right]\qquad\mbox{and}\\
Z_{i,n}&:=&D_{i,n}-D_{i,n+1}+\beta (R_{i,n}-R_{i,n+1})\\
&&+\ 4\won\{n\mbox{ is good for $i$}\}.
\end{eqnarray*}
Therefore, in order to prove  (\ref{chin}) we need to 
bound $f_{i,n}(c_{\ref{alll}})$ on $\{\tau_i+n<\tau_{i+1}\}$ from above away from 1.
By Taylor's expansion, 
\begin{equation} 
f_{i,n}(c_{\ref{alll}})
\leq 1+c_{\ref{alll}} E[Z_{i,n}\ |\ \F_{\tau_i,n}]+
(c_{\ref{alll}}c_{\ref{boun}})^2 \exp(c_{\ref{alll}} c_{\ref{boun}})/2,\label{emo}
\end{equation}
\abel{emo}
where $\con{boun}:=1+\beta+4$ is an upper bound for $Z_{i,n}$.
On $\{\tau_i+n<\tau_{i+1}\},$ due to definition (\ref{bega}) of $\beta$,
\begin{eqnarray*}
E[Z_{i,n}\ |\ \F_{\tau_i,n}]
&=&E[D_{i,n}-D_{i,n+1}+R_{i,n}-R_{i,n+1}\ |\  \F_{\tau_i+n}]\\
&&+\ 4\left(\pi \sqrt{8}
E[R_{i,n}-R_{i,n+1}\mid  \F_{\tau_i+n}]+\won
\{n\mbox{ is good for $i$}\}\right) \\
&\leq&-c_{\ref{ekg}}
\end{eqnarray*}
$P$-a.s. by virtue of Lemma \ref{neg}  (\ref{blitz}), (\ref{jodea}), and definition (\ref{lean}).
Consequently, we may and do choose $c_{\ref{alll}}>0$ small enough such  that
on $\{\tau_i+n<\tau_{i+1}\}$ the right hand side of (\ref{emo}) is $P$-a.s.\ less than a number strictly 
smaller than 1, which we call $e^{-c_{\ref{ofer}}}$, thus showing (\ref{chin}).
By induction over $n$ we obtain
\[
E\left[M_{i,n}\ |\ \F_{\tau_i}\right]\leq \exp(-c_{\ref{ofer}}(n-1))
E\left[M_{i,1}\ |\ \F_{\tau_i}\right]\qquad \mbox{$P$-a.s.}
\]
for all $n\geq 1$. Since $M_{i,1}\leq e^{c_{\ref{alll}} c_{\ref{boun}}}$ this finishes the proof.
One could do better by estimating $M_{i,1}$ more carefully, thus getting rid of the
constant $c_{\ref{elon}}$, but we do not need it.
\end{proof}
%%%%%%%%%%%%%%%%%%%%%%%%%%%%%%%%%%%%%%%%%%%%%%%%%%%%%%%%%%%%%%%%%%%%

\section{Exponential moments of $\tau_{i+1}-\tau_i$.}

We denote the differences between two successive finite ladder points of $(d_n)_n$ by
$\Delta_i:=\tau_{i+1}-\tau_i$ for $i\geq 0$. 

\begin{prop}\label{perv}
$\tau_i<\infty\ P$-a.s.\ for all $i\geq 0$.
Moreover, 
there are constants   $1\leq c_{\ref{elon}}<\infty$ and
$\con{posi}>0$
such that for all $i\geq 0$ and $n\geq 0$, $P$-a.s.,
\begin{eqnarray}
E\left[\exp(c_{\ref{posi}}\Delta_i)\mid \F_{\tau_i}\right]&\leq&c_{\ref{elon}}\label{kleine}\qquad \mbox{and}\\
P[\Delta_i\geq n\ |\ \F_{\tau_i}]&\leq& c_{\ref{elon}}\exp(-c_{\ref{posi}}n).
\label{ur}
\end{eqnarray}
\end{prop}\abel{perv}\abel{posi}\abel{posi}\abel{ur}\abel{kleine}
\begin{proof}
We only need to show that there are  constants  $1\leq c_{\ref{elon}}<\infty$ and  $c_{\ref{posi}}>0$ such that (\ref{kleine})
holds for all $i\geq 1$ with $\tau_i<\infty$. Indeed,
the case $i=0$ is trivial since $\tau_0=0,\ \tau_1=1$ and hence $\Delta_0=1$.
Moreover, (\ref{ur}) follows from (\ref{kleine}) by  Chebyshev's inequality
and since (\ref{kleine}) implies $\Delta_i<\infty$, 
we then have $\tau_i=\Delta_0+\ldots+\Delta_{i-1}<\infty$ as well.

Fix $i\geq 1$. For the proof of (\ref{kleine})
we first show that with  $c_{\ref{alll}}, c_{\ref{ofer}}$ and $c_{\ref{elon}}$ according to 
Lemma \ref{ftau} and  $\con{witri}:=c_{\ref{alll}}(1+2\beta)/c_{\ref{ofer}}$ 
we have $P$-a.s. for all $i\geq 1, n\geq 0$,
\begin{eqnarray}\label{weak}
P\left[\Delta_i>n\ |\ \F_{\tau_i}\right]&\leq& 
c_{\ref{elon}}\exp(c_{\ref{ofer}}(c_{\ref{witri}}d_{\tau_i}-n)),\quad\mbox{and}\\
P\left[\tau_i+n<\tau_{i+1}\wedge \ga_{i+1}\ |\ \F_{\tau_i}\right]&\leq& 
c_{\ref{elon}}\exp(-c_{\ref{ofer}}n).\label{strong}
\end{eqnarray}\abel{weak}\abel{strong}
 We  shall  show later how these auxiliary estimates imply (\ref{kleine}).
For $n=0$, (\ref{weak}) and (\ref{strong}) are true since $c_{\ref{elon}}\geq 1$.
Fix $n\geq 1$.
By (\ref{nochmal}) and the Pythagorean theorem, 
$D_{i,n}\leq R_{i,n}\leq  d_{\tau_i}$ on the event $\{\tau_i+n<\tau_{i+1}\}$.
Hence due to Lemma \ref{ftau} (\ref{chen}) $P$-a.s.\ for all $n\geq 1$,
\[\exp(-c_{\ref{alll}}(1+2\beta) d_{\tau_i})P[\tau_i+n<\tau_{i+1}\ |\ \F_{\tau_i}]\leq
c_{\ref{elon}}\exp(-c_{\ref{ofer}}n),
\]
which is equivalent to (\ref{weak}).
The second auxiliary statement (\ref{strong}) follows from  Lemma \ref{ftau} 
(\ref{chen}) once we have shown that on the event $\{\tau_i+n<\tau_{i+1}\wedge \ga_{i+1}\}$,
\begin{equation}\label{justi}
D_{i,n}+\beta(R_{i,n}-R_{i,0})-4\sum_{j=0}^{n-1}\won\{\mbox{$j$ is good for $i$}\}\leq 0.
\end{equation}\abel{justi}
First we will show that on $\{\tau_i+n<\tau_{i+1}\wedge \ga_{i+1}\}$,
\begin{equation}\label{fair}
D_{i,n}+\beta(R_{i,n}-R_{i,0})\leq 2(D_{i,n}-\|X_{\tau_i}-Z_{i,n}\|).
\end{equation}\abel{fair}
This is done by brute force. For abbreviation we set $d:=d_{\tau_1}/2, y:=D_{i,n}$ and $x:=\|X_{\tau_i}-Z_{i,n}\|$
and note that on  $\{\tau_i+n<\tau_{i+1}\wedge \ga_{i+1}\}$ we have $x,y\in [0,\ga d]$.
Observe that $x$ and $y$ play the role of cartesian coordinates of $X_{\tau_i+n}$, see Figure \ref{isch}.
Using $R_{i,0}=d/2$ and
$R_{i,n}=\sqrt{(d/2-x)^2+y^2}$ we see that (\ref{fair}) is equivalent to
\[\beta\sqrt{(d/2-x)^2+y^2}\leq y-2x+\beta d/2.\]
Both sides of this inequality are nonnegative since $x$ is less than $\ga d$, which is tiny compared to $\beta d$.
Taking the square and rearranging shows that
(\ref{fair}) is equivalent to
\begin{equation}\label{dave}
x(4x-4y-2\beta d-\beta^2 x+\beta^2 d)+y(y+\beta d-\beta^2y)\geq 0.
\end{equation}\abel{dave}
Since $x,y\in [0,\ga d]$ and $\beta\ga=1/2$, see (\ref{bega}), the terms $\beta^2d$ 
in the first bracket and  $\beta d$ in the second bracket
are the dominant terms, respectively,  which shows that (\ref{dave}) and thus (\ref{fair}) holds.
For the proof of (\ref{justi}) it therefore suffices to show that on $\{\tau_i+n<\tau_{i+1}\wedge \ga_{i+1}\}$,
\begin{equation}\label{Honor}
D_{i,n}-\|X_{\tau_i}-Z_{i,n}\|\leq 2\sum_{j=0}^{n-1}\won\{\mbox{$j$ is good for $i$}\}.
\end{equation}\abel{Honor}
Both $D_{i,j}$ and $-\|X_{\tau_i}-Z_{i,j}\|$ can increase by at most 1 if $j$ increases by 1.
Therefore, the left hand side of (\ref{Honor}) is less than or equal to $2\#J$ where
\[J:=\{0\leq j< n\mid \forall 0\leq m<j: D_{i,m}-\|X_{\tau_i}-Z_{i,m}\|\leq D_{i,j}-\|X_{\tau_i}-Z_{i,j}\|
\}.\]
Hence it suffices to show that the elements of $J$ are good for $i$.
Note that $j=0\in J$ is good for $i$ by definition of being good. 
So fix $1\leq j\in J$. 
 By Lemma \ref{neg} (\ref{blitz}) it is enough to show that
$\sin\phi_{1,i,j}\geq 2^{-1/2}$, that is
\begin{equation}\label{coffee}
\phi_{1,i,j}\in [\pi/4, 3\pi/4].
\end{equation}\abel{coffee} 
On the one hand, $\phi_1-\psi_1$ is close to $\pi/2$, as can be seen
in Figure \ref{isch}. 
\begin{figure}[t]
\epsfig{file=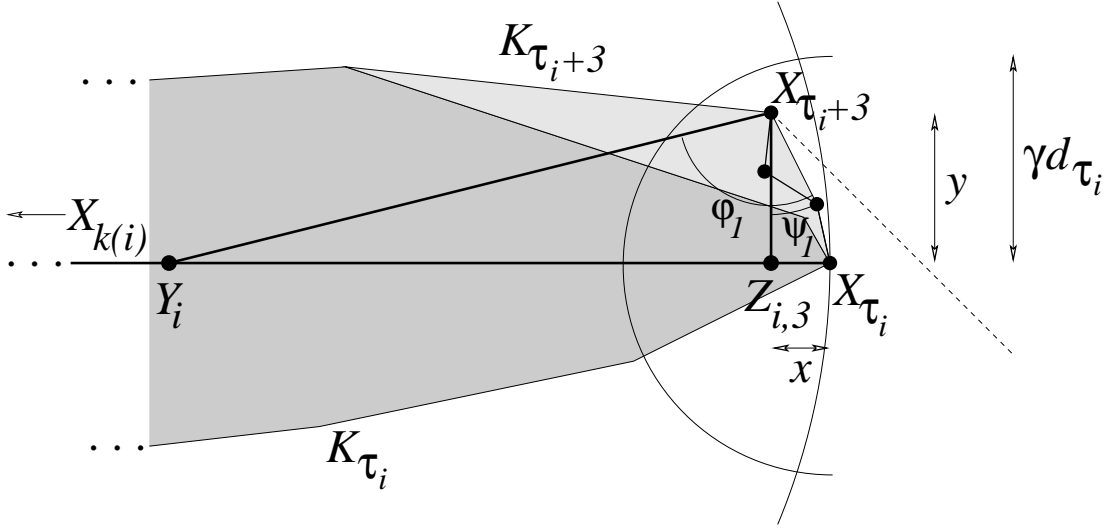,height=7cm,angle=0}
\caption{\foot  
The darkly shaded convex hull $K_{\tau_i}$ at time $\tau_i$  
has been enlarged after three steps by the lightly  shaded part. $j=3$ is good
for $i$ since it satisfies the sufficient criterion $\psi_{1,i,3}\leq \pi/4$, see (\ref{ill}), 
which corresponds to the fact that
the dashed line, which intersects the horizontal axis at an angle of $\pi/4$,
does not intersect $K_{\tau_i+3}^o$. $j=1$ is also good  for $i$ for the same reason, 
while $j=2$ might be good for $i$ but fails to satisfy the sufficient condition (\ref{ill}), since the
corresponding dashed line starting in $X_{\tau_i+2}$ would have intersected $K_{\tau_i+2}^o$.
}
\label{isch}
\end{figure}\abel{isch}
More precisely,
\begin{eqnarray*}
\sin(\phi_1-\psi_1)&=&\frac{\|Y_i-Z_{i,j}\|}{\|Y_i-X_{\tau_i+j}\|}\geq 
\frac{\|Y_i-X_{\tau_i}\|-\|X_{\tau_i}-Z_{i,j}\|}{
\|Y_i-X_{\tau_i}\|+\|X_{\tau_i}-X_{\tau_i+j}\|}\\
 &\geq& 
\frac{d_{\tau_i}/2-\ga d_{\tau_i}}{d_{\tau_i}/2+\ga d_{\tau_i}}\ =\ \frac{1-2\ga}{1+2\ga}\ \geq\ \frac{1}{\sqrt{2}}
\ =\ \sin\frac{\pi}{4}.
\end{eqnarray*}
Since $0\leq \psi_1\leq \phi_1$ this implies $\phi_1\geq \pi/4$.
On the other hand, $\phi_1-\psi_1\leq \pi/2$, see (\ref{fax}). Hence all that remains 
 to be shown for the completion of the
proof of (\ref{coffee}) and (\ref{strong}) is  that 
\begin{equation}\label{ill}
\psi_{1,i,j}\leq \pi/4.
\end{equation}\abel{ill} 
Consider the half line (dashed in Figure \ref{isch})  starting at 
$X_{\tau_i+j}$ which includes an angle of $\pi/4$ with
$\overline{X_{\tau_i+j},Z_{i,j}}$ 
that contains $s_{1,i,j}$.
We claim that this line does not intersect $K_{\tau_i+j}^o$. This would
imply (\ref{ill}).
 To prove this claim observe that  for any $c>0$ the 
set of possible values for $X_{\tau_i+m}$ with $D_{i,m}-\|X_{\tau_i}-Z_{i,m}\|=c$  
is a line 
parallel to the half line just described.
Since $j\in J$
the walker did not cross between time $\tau_i$ and time $\tau_i+j-1$ the dashed line
passing through  $X_{\tau_i+j}$. Consequently, it suffices  to show that the dashed line does not intersect
$K_{\tau_i}^o$. If it did intersect  $K_{\tau_i}^o$ then this 
 would force the walker  on its way  from
$X_{\tau_i}$ to $X_{\tau_i+j}$ to cross the dashed line strictly before time $\tau_i+j$, which is impossible as we just saw.
This completes the proof of  (\ref{strong}).

%%%%%%%%%%%%%%%%%%%%%%%%%%%%%%%%%%%%%%%%%

Finally, we demonstrate how (\ref{weak}) and (\ref{strong}) imply (\ref{kleine}) with 
\begin{equation}\label{flu}
c_{\ref{posi}}:=
\frac{\ga c_{\ref{ofer}}}{c_{\ref{witri}}+\ga}.
\end{equation}
\abel{flu}
We distinguish three cases by partitioning $\Omega$ into
three elements of $\F_{\tau_i}$:
\[
\{n\leq\ga d_{\tau_i}\},\quad 
\{\ga d_{\tau_i}<n<(c_{\ref{witri}}+\ga)d_{\tau_i}\}\quad
\mbox{and}\quad \{(c_{\ref{witri}}+\ga)d_{\tau_i}\leq n\}.
\]
 Note that 
\begin{equation}\label{bombs}
\ga_{i+1}\geq\tau_i+\lceil \ga d_{\tau_i}\rceil
\end{equation}\abel{bombs}
since the walker takes steps of length one. Therefore, on $\{n\leq\ga d_{\tau_i}\}$, 
\begin{eqnarray*}
P\left[\Delta_i>n\ |\ \F_{\tau_i}\right]\ &\stackrel{(\ref{bombs})}{=}&
P\left[\tau_i+n<\tau_{i+1}\wedge \ga_{i+1}\ |\ \F_{\tau_i}\right]
\\
 &\stackrel{(\ref{strong})}{\leq}& c_{\ref{elon}}\exp(-c_{\ref{ofer}}n)\ \stackrel{(\ref{flu})}{\leq}\ 
c_{\ref{elon}}\exp(-c_{\ref{posi}}n).
\end{eqnarray*}
On $\{\ga d_{\tau_i}<n<(c_{\ref{witri}}+\ga)d_{\tau_i}\}$, 
\begin{eqnarray*}\lefteqn{
P\left[\Delta_i>n\ |\ \F_{\tau_i}\right]}\\
&\leq& P\left[\tau_i+\lceil \ga d_{\tau_i}\rceil<\tau_{i+1}\ |\ \F_{\tau_i}\right]\ 
\stackrel{(\ref{bombs})}{=}\ 
P\left[\tau_i+\lceil\ga d_{\tau_i}\rceil<\tau_{i+1}\wedge \ga_{i+1}\ 
|\ \F_{\tau_i}\right]\\
&=&\sum_{k\geq 1}P\left[\tau_i+k<\tau_{i+1}\wedge \ga_{i+1}\ 
|\ \F_{\tau_i}\right]\won\{\lceil\ga d_{\tau_i}\rceil=k\}\\
&\stackrel{(\ref{strong})}{\leq}& c_{\ref{elon}}\exp(-c_{\ref{ofer}}\ga d_{\tau_i})\ 
\stackrel{(\ref{flu})}{\leq}\ c_{\ref{elon}}\exp(-c_{\ref{posi}}n).
\end{eqnarray*}
Finally, on $\{(c_{\ref{witri}}+\ga)d_{\tau_i}\leq n\}$, 
\begin{eqnarray*}
P\left[\Delta_i>n\ |\ \F_{\tau_i}\right]
&\stackrel{(\ref{weak})}{\leq}&c_{\ref{elon}}\exp(c_{\ref{ofer}}(c_{\ref{witri}}d_{\tau_i}-n))\ 
\stackrel{(\ref{flu})}{\leq}\ c_{\ref{elon}}\exp(-c_{\ref{posi}}n),
\end{eqnarray*}
where the last inequality can easily be checked.
\end{proof}
%%%%%%%%%%%%%%%%%%%%%%%%%%%%%%%%%%%%%%%%%%%%%%%%%%%%%%%%%%%%%%%%%%%%%%%%%

\section{Linear growth of the diameter and proof of Theorem \ref{kalt}}
The following result (with $c=0$) implies that $(d_n)_n$ has a positive lim inf speed.
\begin{lemma}\label{ai}
There are constants $\con{lai}>0$ and $\con{map}<\infty$ such that for all $n\geq 0$ and all $c\in[0,1[$,
\begin{equation}\label{aussi}
E[\exp(d_{\lfloor cn\rfloor}-d_n)]\leq c_{\ref{map}}(n+1)\exp(-c_{\ref{lai}}(1-c)n).
\end{equation}\abel{aussi}
\end{lemma}\abel{ai}\abel{lai}\abel{map}
%%%%%%%%%%%%%%%%%%%%%%%%%%%%%%%%%%%%%%%%%%%%%%%%%%%%%%%%%%%%%%%%%%%%%%
For the proof of this lemma and of Theorem \ref{kalt} we need the following definition:
  Given $n\geq 0$ let
$i_n:=\sup\{i\geq 0\mid \tau_i\leq n\}$. Note that 
\begin{equation}\label{bus}
 d_{\tau_{i_n}}=d_n\quad\mbox{and}\quad i_n\leq \tau_{i_n}\leq n<\tau_{i_{n+1}}.
\end{equation}\abel{bus}
\begin{proof}[Proof of Lemma \ref{ai}.]
The case $n=0$ is trivial. Now fix $n\geq 1$, 
$c\in[0,1[$ and set $m=\lfloor cn\rfloor$,
\begin{equation}\label{tas}
\con{pots}:=\frac{1-c}{2}>0\qquad\mbox{and}\qquad 
\con{porta}:=\frac{c_{\ref{posi}}c_{\ref{pots}}}{2\ln c_{\ref{elon}}}>0,
\end{equation} \abel{tas}
where  $c_{\ref{elon}}$ and $c_{\ref{posi}}$ are according to Proposition \ref{perv}.
A simple union bound yields
\begin{eqnarray*}
E[e^{d_m-d_n}]&\leq& {\rm I}+{\rm II}+{\rm III},\qquad\mbox{where}\\
{\rm I}&:=&P[\tau_{i_m+1}-m\geq c_{\ref{pots}}n],\\
{\rm II}&:=&P[\tau_{i_m+1}-m< c_{\ref{pots}}n,\ i_n<i_m+\lceil c_{\ref{porta}}n\rceil]\qquad\mbox{and}\\
{\rm III}&:=& E\left[\exp\left(d_{m}-d_{n}
\right),\ i_n\geq i_m+\lceil c_{\ref{porta}}n\rceil\right],
\end{eqnarray*}
see also Figure \ref{line}.
%%%%%%%%%%%%%%%%%%%%%%%%%%%%%%%%%%%%%%%%%%%%%%%%%%%%%%%%%%%%%%%%%%%%%%
\begin{figure}[t]
\epsfig{file=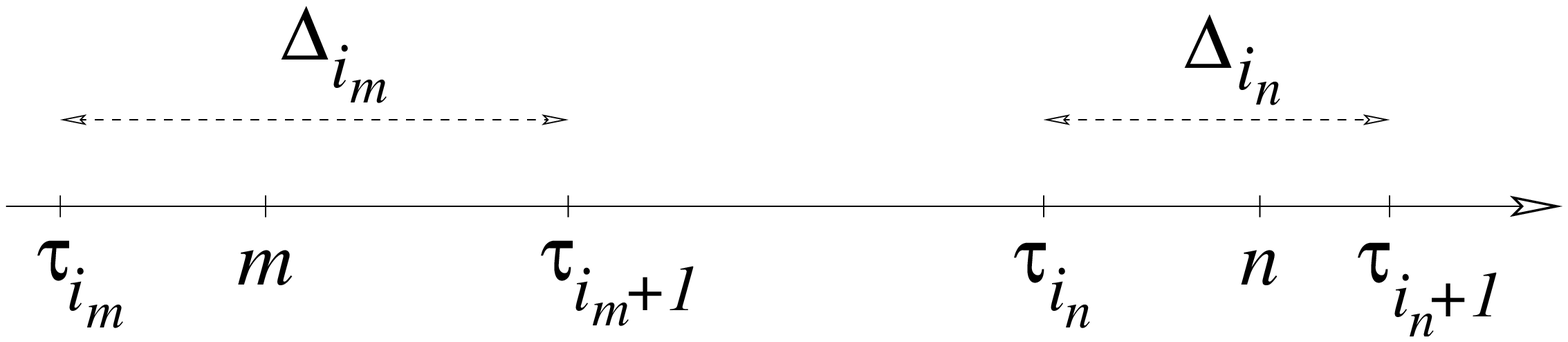,height=2cm,angle=0}
\caption{\foot  
}\label{line}
\end{figure}\abel{line}
%%%%%%%%%%%%%%%%%%%%%%%%%%%%%%%%%%%%%%%%%%%%%%%%%%%%%%%%%%%%%%%%%%%%%%
Here term I corresponds to the situation in which 
after time $m$ the diameter does not increase for an untypical long while.
 Term II handles the case in which the diameter does increase
shortly after time $m$, as it should,  but
not often enough in
the remaining time until $n$.
The third term III considers the original random variable on the typical event that the 
number of times at which the
diameter increases
is at least  proportional to $n$ with a constant of proportionality not too small. 

It suffices  to show that each of these three terms decays as $n\to\infty$ in the way stated in (\ref{aussi})
with constants $c_{\ref{lai}}$ and $c_{\ref{map}}$ independent of $c$.
As for the first term, 
\[
{\rm I}\leq P[\Delta_{i_m}\geq \lceil c_{\ref{pots}} n\rceil]
\stackrel{(\ref{bus})}{\leq}\sum_{i=0}^{m}P[i_m=i,\ \Delta_{i}\geq \lceil c_{\ref{pots}} n\rceil]\\
\stackrel{(\ref{ur})}{\leq} c_{\ref{elon}}(n+1)e^{-c_{\ref{posi}}c_{\ref{pots}} n},
\]
which is an upper bound like the one requested in (\ref{aussi}).
The second term is estimated as follows.
\begin{eqnarray}
{\rm II}&\stackrel{(\ref{bus})}{=}& P\left[\tau_{i_m+1}<\lfloor cn\rfloor +c_{\ref{pots}}n,\ n<\tau_{i_n+1}\leq 
\tau_{i_m+\lceil c_{\ref{porta}}n\rceil}
\right]\nonumber\\ 
&\stackrel{(\ref{tas})}{\leq}& P\left[ \tau_{i_m+\lceil c_{\ref{porta}}n\rceil}-\tau_{i_m+1}\geq 
(1-c-c_{\ref{pots}})n=c_{\ref{pots}}n\right]\nonumber\\
&\leq&  E\left[\exp\left(c_{\ref{posi}}\left(
 \tau_{i_m+\lceil c_{\ref{porta}}n\rceil}-\tau_{i_m+1}-c_{\ref{pots}}n\right)\right)\right]\nonumber\\
&=&e^{-c_{\ref{posi}}c_{\ref{pots}}n}\sum_{k\geq 1} E\left[\exp\left(c_{\ref{posi}}\left(
\tau_{k+\lceil c_{\ref{porta}}n\rceil-1}-\tau_{k}\right)
\right), i_m+1=k\right]\nonumber\\
&=&e^{-c_{\ref{posi}}c_{\ref{pots}}n}\sum_{k\geq 1} E\left[\prod_{i=0}^{\lceil c_{\ref{porta}}n\rceil-2}
\exp\left(c_{\ref{posi}}\Delta_{k+i}
\right), i_m+1=k\right].\label{wo}
\end{eqnarray}\abel{wo}
Note that $\{i_m+1=k\}$ is the event that $\tau_k$ is the first time after time $m$ at which 
the diameter increases. Therefore, 
\begin{equation}\label{bank}
\{i_m+1=k\}\in\F_{\tau_k}.
\end{equation}\abel{bank}
Moreover, the increments $\Delta_{k+i}$
are measurable with respect to  $\F_{\tau_{k+i+1}}$. 
Consequently, by conditioning in (\ref{wo}) on  $\F_{\tau_{k+\lceil c_{\ref{porta}}n\rceil-2}}$ and applying
 Proposition \ref{perv} (\ref{kleine}) with $i=k+\lceil c_{\ref{porta}}n\rceil-2$ we conclude 
\[{\rm II}\leq e^{-c_{\ref{posi}}c_{\ref{pots}}n}c_{\ref{elon}}
\sum_{k\geq 1} E\left[\prod_{i=0}^{\lceil c_{\ref{porta}}n\rceil-3}
\exp\left(c_{\ref{posi}}\left(\Delta_{k+i}\right)
\right), i_m+1=k\right].\]
Continuing in this way we obtain by induction after $\lceil c_{\ref{porta}}n\rceil-1$ steps, 
\[{\rm II}\leq e^{-c_{\ref{posi}}c_{\ref{pots}}n}c_{\ref{elon}}^{\lceil c_{\ref{porta}}n\rceil-1}
\ \leq\ e^{-c_{\ref{posi}}c_{\ref{pots}}n}c_{\ref{elon}}^{ c_{\ref{porta}}n}
\ \stackrel{(\ref{tas})}{=}\  e^{-(c_{\ref{posi}}/4)(1-c)n},\]
 which is again of the form required in (\ref{aussi}). 

In order demonstrate that also the third term III behaves properly we will 
show that the increments $d_{\tau_{i+1}}-d_{\tau_i}, i\geq 1,$ 
have a uniformly positive chance of being larger than a fixed constant, say 1/2,
independently of the past.
More precisely, we  may assume that the process $(X_n)_n$ is generated in the following way:
There are i.i.d.\ random variables $U_{n,k},\ n\geq 0, k\geq 0,$ uniformly distributed
on the unit circle centered in 0 such that 
$X_{n+1}=X_n+U_{n,k},$
where $k$ is the smallest integer such that 
$\overline{X_n, X_n+U_{n,k}}$
does not intersect $K_n^o$.
Then for any $i\geq 1$, 
by definition of $\tau_i$,
\begin{eqnarray}\nonumber
\{d_{\tau_{i+1}}\geq d_{\tau_i}+1/2\}&\supseteq& \{d_{\tau_{i}+1}\geq d_{\tau_i}+1/2\}\\
&\supseteq& \left\{U_{\tau_i,0}\cdot \frac{X_{\tau_i}-X_{k(i)}}{d_{\tau_i}}
\geq \frac{1}{2}\right\}=:A_i,\label{food}
\end{eqnarray}\abel{food}
see Figure \ref{half}.
%%%%%%%%%%%%%%%%%%%%%%%%%%%%%%%%%%%%%%%%%%%%%%%%%%%%%%%%%%%%%%
\begin{figure}[t]
\epsfig{file=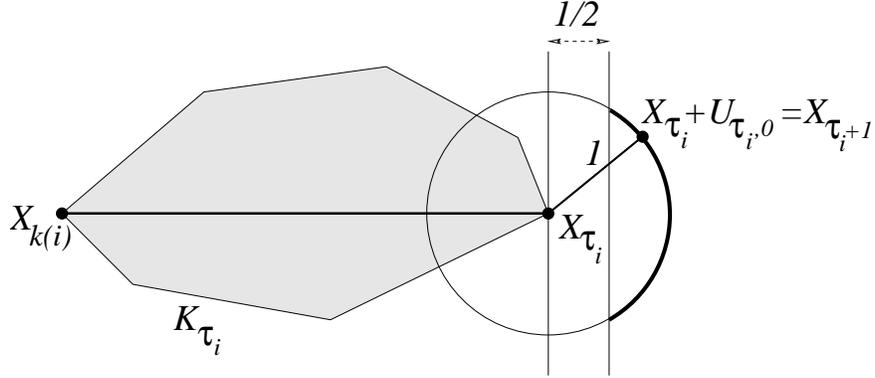,height=5cm,angle=0}
\caption{\foot  The event $A_i$ occurs if the first trial point sampled lies on the bold arc.
}\label{half}
\end{figure}\abel{half}
%%%%%%%%%%%%%%%%%%%%%%%%%%%%%%%%%%%%%%%%%%%%%%%%%%%%%%%%%%%%%%
Here (\ref{food}) holds for the following reason: 
Observe that for all $x\in K_{\tau_i}$,
\begin{eqnarray*}
(X_{\tau_i}-X_{k(i)})\cdot\frac{X_{\tau_i}-X_{k(i)}}{\|X_{\tau_i}-X_{k(i)}\|}&=&d_{\tau_i}\ \geq\
\|x-X_{k(i)}\|\\
& \geq& 
(x-X_{k(i)})\cdot\frac{X_{\tau_i}-X_{k(i)}}{\|X_{\tau_i}-X_{k(i)}\|}
\end{eqnarray*}
by Cauchy-Schwarz and thus 
\begin{equation}\label{show}
x\cdot (X_{\tau_i}-X_{k(i)})\leq X_{\tau_i}\cdot (X_{\tau_i}-X_{k(i)})\qquad (x\in K_{\tau_i}).
\end{equation}\abel{show}
However, on $A_i$, 
\[
(X_{\tau_i}+U_{\tau_i,0})\cdot (X_{\tau_i}-X_{k(i)})\geq 
X_{\tau_i}\cdot (X_{\tau_i}-X_{k(i)})+d_{\tau_i}/2> X_{\tau_i}\cdot (X_{\tau_i}-X_{k(i)}).
\]
Along with (\ref{show}) this shows that 
 $X_{\tau_i}+U_{\tau_i,0}$ and $K_{\tau_i}$
are lying  on opposite sides of the line passing orthogonally to $\overline{X_{\tau_i},X_{k(i)}}$
through  $X_{\tau_i}$. Therefore, 
$\overline{X_{\tau_i},X_{\tau_i}+U_{\tau_i,0}}$
does not intersect $K_{\tau_i}^o$.  Consequently, on $A_i$,
$X_{\tau_i+1}=X_{\tau_i}+U_{\tau_i,0}$ and thus by Cauchy-Schwarz
\begin{eqnarray*}
d_{\tau_i+1}&\geq& \|X_{\tau_i+1}-X_{k(i)}\|\ \geq\
 (X_{\tau_i}+U_{\tau_i,0}-X_{k(i)})\cdot \frac{X_{\tau_i}-X_{k(i)}}
{d_{\tau_i}}
\ \geq\ d_{\tau_i}+\frac{1}{2},
\end{eqnarray*}
which completes the proof of (\ref{food}).
Therefore, for all $1\leq j_1\leq j_2$,
\begin{eqnarray}
d_{\tau_{j_2}}-d_{\tau_{j_1}}&\geq& \frac{1}{2}\sum_{i=j_1}^{j_2-1}\won\{d_{\tau_{i+1}}\geq 
d_{\tau_i}+1/2\}
\ \geq\  
 \frac{1}{2}\sum_{i=j_1}^{j_2-1}\won\{A_{i}\}.\label{sars}
\end{eqnarray}\abel{sars}
This estimate will be useful since the 
 random variables 
\begin{equation}\label{water}
\won\{A_i\}\ (i\geq 1)\quad
\mbox{are i.i.d.\ with $P[A_i]>0$}.
\end{equation}\abel{water}
Indeed, let $\widetilde{\F}_n\ (n\geq 0)$ be the $\sigma$-field generated by
$U_{m,k},\ 0\leq m<n, 0\leq k$. Because of $\F_n\subseteq \widetilde{\F}_n$ we have 
$A_j\in  \widetilde{\F}_{\tau_i}$ for all $1\leq j<i$. 
Moreover, since the uniform distribution on the unit circle is invariant under
rotations, 
\begin{equation}\label{nochnlabel}
\mbox{$A_i$ is independent of  $\widetilde{\F}_{\tau_i}\ (i\geq 1)$}
\end{equation}\abel{nochnlabel}
 and $P[A_i\mid \widetilde{\F}_{\tau_i}]=P[A_i]$
is just the length of the bold circle segment shown in Figure \ref{half} divided by $2\pi$. This implies
(\ref{water}).
Now we  estimate III by
\begin{eqnarray}\nonumber
{\rm III}&\stackrel{(\ref{bus})}{=}&
 E\left[\exp\left(d_{\tau_{i_m}}-d_{\tau_{i_n}}
\right),\ i_n\geq i_m+\lceil c_{\ref{porta}}n\rceil\right]\\ \nonumber
&\leq&  E\left[\exp\left(d_{\tau_{i_m+1}}-d_{\tau_{i_n}}
\right),\ i_n\geq (i_m+1)+ \lceil c_{\ref{porta}}n\rceil-1\right]\\ \nonumber
&\leq&  \sum_{k\geq 1} E\left[\exp\left(d_{\tau_{k}}-d_{\tau_{k+ \lceil c_{\ref{porta}}n\rceil-1}}
\right),\ i_m+1=k\right]\\ \label{right}
&\stackrel{(\ref{sars})}{\leq}& \sum_{k\geq 1}E\left[\exp\left(-\frac{1}{2}\sum_{i=k}^{k+ 
\lceil c_{\ref{porta}}n\rceil-2}
\won\{A_i\}\right), i_m+1=k\right].
\end{eqnarray}\abel{right}
As seen in (\ref{bank}),  $\{i_m+1=k\}\in \F_{\tau_k}\subseteq \widetilde{\F}_{\tau_k}$.
Therefore, after conditioning in (\ref{right}) on $\F_{\tau_k}$, we see with the help of (\ref{nochnlabel})
for $i\geq k$ and 
(\ref{water})  that the right hand side 
of (\ref{right}) equals
\[ E\left[\exp\left(-\frac{1}{2}A_1\right)\right]^{\lceil c_{\ref{porta}}n\rceil-1},\]
which  decays as required in (\ref{aussi}), see (\ref{tas}).     
\end{proof}

%%%%%%%%%%%%%%%%%%%%%%%%%%%%%%%%%%%%%%%%%%%%%%%%%%%%%%%%%%%%%%%%%%%%%%%%%%

Lemma \ref{ai} directly implies a weaker version of Theorem 
\ref{kalt} in which  $\|X_n\|$ is replaced by $d_n$. For the full statement we need the following additional argument.
\begin{proof}[Proof of Theorem \ref{kalt}]
(\ref{part}) follows from (\ref{kalt}) by the Borel-Cantelli lemma.
For the proof of (\ref{kalt})
pick  $c_{\ref{lai}}$ and $c_{\ref{map}}$ according to Lemma \ref{ai} and 
choose $\con{wer}>0$ and $c_{\ref{one}}>0$ small enough such that
\begin{equation}\label{korea}
2c_{\ref{wer}}-c_{\ref{lai}}<0\qquad\mbox{and}\qquad 2c_{\ref{one}}-c_{\ref{lai}}(c_{\ref{wer}}-c_{\ref{one}})<0.
\end{equation}\abel{korea}
We denote by
$M_n:=\max\{\|X_m\|\ |\ m\leq n\}$ the walker's maximal distance from the origin 
by time $n$. Note that $M_n$ and $d_n$ are related via
\begin{equation}\label{mob}
M_n\leq d_n\leq 2 M_n\quad\mbox{for all $n\geq 0$}
\end{equation}\abel{mob}
because of $X_0=0$.
 By a union bound for any $n\geq 0$,
\begin{eqnarray}
P[\|X_n\|\leq c_{\ref{one}} n] 
&\leq& P[\Delta_{i_n}\geq c_{\ref{one}} n]+P[M_n\leq c_{\ref{wer}}n]+
P[B_n],\ \mbox{where}\label{andrew}
\\
B_n&:=& \{\Delta_{i_n}<c_{\ref{one}} n,\ M_n>c_{\ref{wer}}n,\ \|X_n\|\leq c_{\ref{one}} n\}.
\nonumber
\end{eqnarray}\abel{andrew}
It suffices to show that each one of the three terms on the right hand side
of (\ref{andrew}) decays exponentially fast in $n$. As for the first term, 
\[P[\Delta_{i_n}\geq c_{\ref{one}} n]\stackrel{(\ref{bus})}{\leq}\sum_{i=0}^{n}
P[i_n=i, \Delta_{i}\geq c_{\ref{one}} n]
\stackrel{(\ref{ur})}{\leq} c_{\ref{elon}}(n+1)e^{-c_{\ref{posi}}c_{\ref{one}} n},\]
which decays exponentially fast in $n$ indeed.  So does the second term in (\ref{andrew}) since
by Chebyshev's inequality,
\[P[M_n\leq c_{\ref{wer}}n]\stackrel{(\ref{mob})}{\leq} P[d_n\leq 2c_{\ref{wer}}n]\leq
e^{2c_{\ref{wer}}n}E[e^{-d_n}]\stackrel{(\ref{aussi})}{\leq} 
c_{\ref{map}}(n+1)e^{(2c_{\ref{wer}}-c_{\ref{lai}})n},\]
which decays exponentially fast due to the choice of $c_{\ref{wer}}$ in (\ref{korea}).
Finally, we are going to bound the third term  in (\ref{andrew}), $P[B_n]$. 
Define the ladder times $(\mu_j)_j$ of the process $(M_n)_{n\geq 0}$ recursively by
\[\mu_0:=0\qquad\mbox{and}\qquad\mu_{j+1}:=\inf\{n>\mu_j\ |\ M_n>M_{\mu_j}\}.\]
In analogy to $(i_n)_n$ for $(\tau_i)_i$ we define for $(\mu_j)_j$ 
the increasing sequence $(j_n)_n$ by $j_n:=\sup\{j\geq 0\mid\mu_j\leq n\}$ and note that
$\mu_{j_n}\leq n< \mu_{j_n+1}$ and $M_n=\|X_{\mu_{j_n}}\|$.
Hence on the event $B_n$,
\[\|X_{\mu_{j_n}}-X_n\|\geq \|X_{\mu_{j_n}}\|-\|X_n\|=M_n-\|X_n\|\geq
(c_{\ref{wer}}-c_{\ref{one}}) n.
\]
Since the walker takes steps of length one, this implies $n-\mu_{j_n}\geq (c_{\ref{wer}}-c_{\ref{one}}) n$
 and therefore, on the event $B_n$,
\begin{equation}\label{imp}
\mu_{j_n}\leq \lfloor (1-c_{\ref{wer}}+c_{\ref{one}}) n\rfloor. 
\end{equation}\abel{imp}
On the other hand, on $B_n$,
\begin{eqnarray*}
d_n&\stackrel{(\ref{bus})}{=}&d_{\tau_{i_n}}\ =\ 
\|X_{\tau_{i_n}}-X_{k(i_n)}\|\ \leq\ \|X_n\|+\|X_n-X_{\tau_{i_n}}\|+\|X_{k(i_n)}\|\\
&\leq&c_{\ref{one}} n+\Delta_{i_n}+M_n\ \leq\ c_{\ref{one}} n+c_{\ref{one}} n+M_{\mu_{j_n}}\\
&\stackrel{(\ref{mob})}{\leq}& 2c_{\ref{one}} n+d_{\mu_{j_n}}
\ \stackrel{(\ref{imp})}{\leq}\ 2c_{\ref{one}} n+d_{
\lfloor (1-c_{\ref{wer}}+c_{\ref{one}}) n\rfloor},
\end{eqnarray*}
where we used in the second inequality again the fact that the steps have length one.
Therefore, by  Chebyshev's inequality and (\ref{aussi}),
\[
P[B_n]\leq P[d_n-d_{\lfloor (1-c_{\ref{wer}}+c_{\ref{one}}) n\rfloor}\leq 2c_{\ref{one}} n] 
\leq c_{\ref{map}}(n+1)e^{(2c_{\ref{one}} -c_{\ref{lai}}(c_{\ref{wer}}-c_{\ref{one}})) n},
\]
which decays exponentially in $n$ due to the choice of $c_{\ref{wer}}$ and $c_{\ref{one}}$ in (\ref{korea}).
\end{proof}

%%%%%%%%%%%%%%%%%%%%%%%%%%%%%%%%%%%%%%%%%%%%%%%%%%%%%%%%%%%%%%%%%%%%%%%%%%

\bibliographystyle{amsalpha}
\vspace*{5mm}
{\sc \small
Department of Mathematics\\
Stanford University\\
Stanford, CA 94305, U.S.A.\\
E-Mail: {\rm zerner@stanford.edu} }
\end{document}